\documentclass[letterpaper, 10 pt, conference]{ieeeconf}  

\usepackage{amsmath}
\usepackage{amssymb}
\usepackage{hyperref}
\usepackage{tikz}
\usepackage{pgfplots}

\IEEEoverridecommandlockouts                              

\title{\LARGE \bf
Generating a robustly stabilizable class of nonlinear systems for a converse optimality problem
}

\author{Rania Tafat$^1$, Thomas Göhrt$^1$ and Stefan Streif$^1$ 
\thanks{$^{1}$ The authors are with the Technische Universit\"at Chemnitz, 09126 Chemnitz, Germany, Automatic Control and System Dynamics Lab; E-mail:
        {\tt\small \{rania.tafat, thomas.goehrt, stefan.streif\}@etit.tu-chemnitz.de}}%
}

\usepackage{graphicx}
\usepackage{caption}
\usepackage{mwe}

\usepackage{xspace}
\usepackage{soul}
\usepackage{mathtools}
\usepackage{makecell}
\usepackage{cite}
\usepackage{nicefrac}
\usepackage{float}
\usepackage{wrapfig}
\usepackage{color}

\newcommand{\ra}{\rightarrow}

\newtheorem{dfn}{Definition}
\newtheorem{asm}{Assumption}
\newtheorem{prp}{Proposition}
\newtheorem{lem}{Lemma}
\newtheorem{thm}{Theorem}
\newtheorem{crl}{Corollary}
\newtheorem{rem}{Remark}

\newcommand{\R}{\ensuremath{\mathbb{R}}}

\newcommand{\X}{\ensuremath{\mathbb{X}}}

\newcommand{\K}{\ensuremath{\mathcal{K}}\xspace}		
\newcommand{\Kinf}{\ensuremath{\mathcal{K}_{\infty}}\xspace}		

\newcommand{\red}[1]{\textcolor{red}{#1}}

\makeatletter
\newcommand{\subalign}[1]{%
	\vcenter{%
		\Let@ \restore@math@cr \default@tag
		\baselineskip\fontdimen10 \scriptfont\tw@
		\advance\baselineskip\fontdimen12 \scriptfont\tw@
		\lineskip\thr@@\fontdimen8 \scriptfont\thr@@
		\lineskiplimit\lineskip
		\ialign{\hfil$\m@th\scriptstyle##$&$\m@th\scriptstyle{}##$\crcr
			#1\crcr
		}%
	}
}

\begin{document}

\maketitle 
\thispagestyle{empty}
\pagestyle{empty}

\begin{abstract}
Converse optimality theory addresses an optimal control problem conversely where the system is unknown and the value function is chosen. 
Previous work treated this problem both in continuous and discrete time and non-extensively considered disturbances.

In this paper, the converse optimality theory is extended to the class of affine systems with disturbances in continuous time while considering norm constraints on both control inputs and disturbances. 
The admissibility theorem and the design of the internal dynamics model are generalized in this context. 

A robust stabilizability condition is added for the initial converse optimality problem using inverse optimality's tool: the robust control Lyapunov function.
A design for nonlinear class of systems that are both robustly stabilizable and globally asymptotically stable in open loop is obtained.
A case study illustrates the presented theory.

\end{abstract}

\section{Introduction}

Consider the following continuous time system:
 \begin{equation}
\dot{x}=f(x)+g_{1}(x) w+g_{2}(x) u,
\label{eq:sys}
\end{equation}
where $x \in \R^n$ are the states, $w \in \R^m$ are the disturbances and $u \in \R^p$ are the inputs. 
The function $f : \R^n \ra \R^n$ with $f(0)=0$ is denoted as internal dynamics, $g_1: \R^n \ra \R^n \times \R^m$ is denoted as disturbance coupling function and $g_2:   \R^n \ra \R^n \times \R^p$ is denoted as input coupling function, respectively. 

Further, consider the following optimization problem \cite{Primbs1996nonlinear}: 
\begin{equation}
\begin{split}
&\inf_{u} \sup _{w} J(x,u, w) \\
J(x,u, w) &=\int_{0}^{\infty}\left(q(x) +u^{\top} u-w^{\top} w\right) d t,
\end{split}
\label{eq:perform}
\end{equation}
where $q: \R^n \ra \R_{\geq 0}$ with $q(0) = 0$ is the utility function.

If disturbances are discarded, this optimal control problem can be calculated after solving a partial differential equation (PDE) known as the Hamilton-Jacobi-Bellman (HJB) equation whose solution is the value function \cite{Bellman1957}.
In the presence of disturbances, a generalized version of the HJB equation that includes a double optimization, is required. This PDE is called the Hamilton-Jacobi-Isaacs (HJI) equation \cite{Isaacs1965}.

Converse HJB or converse optimality was developed in \cite{nevistic1996optimality} as an approach that addresses another aspect of the optimal control problem where the value function is assumed known and the internal dynamics of the system assumed unknown.
A converse optimality problem addresses the question: "Which class of nonlinear systems satisfies the optimal control problem under a chosen value function?"

The Co-HJB approach was first developed to design benchmark class of systems to test control methods \cite{Primbs1996nonlinear} and then  constraints on states and control inputs were considered in  \cite{nevistic1996constrained}. 
All work on converse optimality was developed in the continuous time framework until \cite{Gohrt2020discrete}, where discrete time is addressed and theorems are developed around it.
To the best knowledge of the authors, work on converse optimality and the examples developed in its context only briefly considered disturbances in the framework but assumed free noise case for co-HJB examples \cite{Primbs1996nonlinear}. 

Another perspective for dealing with optimal control problems is inverse optimality where links between stabilizability and optimality are explored for linear systems \cite{cai2017inverse}, \cite{Kalman1964} and nonlinear systems \cite{krstic1999inverse}, \cite{moylan1973nonlinear}, \cite{pauwels2016linear}, \cite{sepulchre2012constructive}.  
In \cite{FreemanBirkhauser1996}, the authors address a framework of constrained control inputs and disturbances+ to extend the control Lyapunov function \cite{ARTSTEIN19831163} for systems with disturbances, which is known as the robust control Lyapunov function (RCLF). 
It is shown that every RCLF solves the HJI PDE associated with a meaningful game and a method to calculate the optimal control using the RCLF, without solving the HJI equation is provided \cite{FreemanSIAM1996}.  

A fundamental difference between the converse and inverse optimality theories stems in their respective frameworks. 
In fact, converse optimality's main result requires unconstrained control inputs and disturbances for exploitation.
In contrast, it is a requirement to consider constraints for control inputs and disturbances for inverse optimal robust stabilization. 
Both theories can be combined by extending the converse optimality theory to the disturbed case in a local approach. This is achieved by setting up an assumption on the state space that allows the use of the RCLF. 
This work aims at finding a systematic approach to construct a class of robust nonlinear systems that ensures the satisfaction of the optimal control problem's condition locally.
In addition, the constructed systems can be used as benchmark example systems to test existing control methods like it is discussed in \cite{Lin2000}, \cite{Peterson2014}, \cite{Wie1992}.

The paper is organized as follows. Section II presents preliminary notations, definitions and assumptions. Section III sets up the robust converse problem and presents the state space assumption. Section IV extends converse optimality for the disturbed case and provides a class of nonlinear systems that solve the robust converse problem. Design and numerical examples are given in Section V. The last section concludes the results.

\section{Preliminaries}

The Euclidean norm is denoted $\left \Vert \cdot \right \Vert : \R^n \ra \R_{\geq 0}$.
The gradient  $\frac{\partial V(x)}{\partial x}$ of a multivariable function $V$ is denoted $V_x$.
The admissible states, control inputs and disturbances are denoted $\X$, $\mathbb{U}$, $\mathbb{W}$ respectively.
The set-valued maps given for control and disturbance constraints are $\mathcal{U}: \X \ra \mathbb{U}$, $\mathcal{W}: \X \ra \mathbb{W}$, respectively. 
The class of continuously differentiable functions is denoted $C^1$. 
\begin{dfn}[ $\mathcal{K}$-function  \cite{FreemanSIAM1996}]
A function $\chi: \R_{\geq 0} \rightarrow \R_{\geq 0}$ is of class $\mathcal{K}$ if $\chi$ is continuous, strictly increasing, and $\chi(0)=0 .$ 
\end{dfn}
\begin{dfn}[ $\mathcal{K}_{\infty}$-function  \cite{FreemanSIAM1996}]
A function $\chi: \R_{\geq 0} \rightarrow \R_{\geq 0}$ is of class $\mathcal{K}_{\infty}$ if $\chi$ is of class $\mathcal{K}$ and satisfies $\chi(r) \rightarrow \infty$ as $r \rightarrow \infty$. 
\end{dfn}
\begin{dfn}[ $\mathcal{KL}$-function  \cite{FreemanSIAM1996}]
A function $\beta: \R_{\geq 0} \times \R_{\geq 0} \rightarrow \R_{\geq 0}$ is of class $\mathcal{K} \mathcal{L}$ if $\beta(\cdot, t)$ is of class $\mathcal{K}$ for each fixed $t \in  \R_{\geq 0}$,  $\beta(r, t)$ decreases to zero as $t \rightarrow \infty$ for each fixed $r \in \R$.
\end{dfn}
\begin{dfn}
The class of continuous positive definite functions $\mathcal{A}(\X)$ on $\X,$ comprises continuous functions $\alpha: \X \ra \R_{\geq 0}$ such that $\alpha(0)=0$ and $\alpha(x)>0$ for $x \in \X \setminus\{0\}$. 
\end{dfn}
\begin{dfn}
The set $\mathcal{A}_{\infty}(\X)$ comprises those functions $\alpha \in \mathcal{A}(\X)$ for which there exists a class $\Kinf$  function $\psi$ such that $\alpha(x) \geq \psi(\Vert x \Vert)$ for all $x \in \X$.
\end{dfn}
\begin{dfn}[Lie derivative \cite{LieDerivative}]
The Lie derivative $L_zV$ of a function $V \in C^{1}$ with respect to a vector field $z$ is denoted $L_{z}V= V_x z$.
\end{dfn}
\begin{dfn}[Saturation operator]
The saturation operator is denoted $\operatorname{sat}_\alpha : \R^n \ra \R^n$ where:
\begin{equation}
\operatorname{sat}_{\alpha}(y)=\left\{\begin{array}{lr}
\alpha \frac{y}{ \Vert y \Vert} & \Vert y \Vert \geq \alpha \\
y &  \Vert y \Vert < \alpha \\
\end{array}\right.
\end{equation}
\end{dfn}

\begin{asm}[\hspace{-0.02em}\cite{FreemanBirkhauser1996}] The set-valued map for control inputs $\mathcal{U}(x)$ is lower semi-continuous with nonempty closed convex values.
\label{asm:U}
 \end{asm}
\begin{asm}[\hspace{-0.02em}\cite{FreemanBirkhauser1996}] The set-valued map for disturbances $w$ 
$\mathcal{W}(x)$ is upper semi-continuous with nonempty compact values.
\label{asm:W}
\end{asm}

\begin{dfn}[RGUAS \cite{FreemanBirkhauser1996}]
\label{dfn:RGUAS}
 Let $\Omega$ be a compact subset of $\X$ such that $0 \in \Omega$. 
The solutions of (\ref{eq:sys}) are robustly globally uniformly asymptotically stable with respect to $\Omega$ (RGUAS-$\Omega$) if there exists a class $\K\mathcal{L}$ function $\beta$ such that for any initial condition $x_{0} \in \mathbb{X}$ and any admissible disturbance, all solutions $x(t)$ starting from $x_{0}$ exist for all $t \geq t_0$ and satisfy $\vert x(t)\vert_{\Omega} \leq  \beta\left(\left|x_{0}\right|_{\Omega}, t\right)$ for all $t \geq 0$. 
\end{dfn}
\begin{dfn}[Robustly stabilizable system \cite{FreemanBirkhauser1996}]
\label{dfn:RSSystems}
A system (\ref{eq:sys}) is robustly stabilizable (RS) if there exists a control law $k$ and a compact set $\Omega \subset \X$ satisfying $0 \in \Omega$ such that the solutions of (\ref{eq:sys}) are RGUAS- $\Omega$.
\end{dfn}
\begin{dfn}[Robust control Lyapunov function \cite{FreemanSIAM1996}]
\label{dfn:RCLF}
A $C^{1}$ function $V \in \mathcal{A}_{\infty}({ \X}) $ is a \textbf{ robust control Lyapunov function} (RCLF) for the system $( \ref{eq:sys})$ if there exists $c \in \R_{\geq 0}$ such that
\begin{equation}
\inf _{u \in \mathcal{U}(x)} \sup _{w \in \mathcal{W}(x)} L_{z} V(x, u, w)<0,
\label{eq:rclf}
\end{equation}
for all $x \in \{ x \in \X: V(x) > c \} .$
\end{dfn}
\begin{prp}[\hspace{-0.02em}\cite{FreemanSIAM1996}]
\label{prp:RCLFexistence}
 If a system (\ref{eq:sys}) satisfies the assumptions \ref{asm:U} and \ref{asm:W} and if there exists a robust control Lyapunov function $V$ for (\ref{eq:sys}), then the system is robustly stabilizable.
\end{prp}
\begin{prp}[\hspace{-0.02em}\cite{FreemanSIAM1996}]
\label{prp:RCLF_HJI}
 Every robust control Lyapunov function solves the HJI equation associated with a meaningful game.
\end{prp}
These propositions combined together allow the use of a RCLF in the locally robust converse optimality context.
\section{Problem Setup }
Consider the nonlinear system that is affine with respect to both control and disturbance such as in (\ref{eq:sys}) along with the following performance objective
\begin{equation}
\begin{split}
&\inf_{u \in \mathcal{U}(x)} \sup_{w \in \mathcal{W}(x)} J(u, w) \\
J(u, w) &=\int_{0}^{\infty}\left(q(x) +u^{\top} u-w^{\top} w\right) d t.
\end{split}
\label{eq:perform}
\end{equation}
The associated HJI equation is:
\begin{equation}
\label{eq:HJI}
\begin{split}
\min_{u \in \mathcal{U}(x)} \max_{w \in \mathcal{W}(x)} \{  V_x^{\top} \left( f(x) + g_1(x)w + g_2(x)u \right) \\ + q(x) + u^\top u - w^\top w \}  = 0,
\end{split}
\end{equation}
where $V: \R^{n} \ra \R_{\geq 0}$ is called the value function.
\begin{asm}
\label{asm:valuefunc}
The value function $V$ is $C^{1}$, in $\mathcal{A}_{\infty}(\X)$ and its gradient satisfies:
\begin{equation}
 \Vert V_x \Vert > 0, \; \forall x \in \X \setminus \{ 0 \}.
\end{equation}
\end{asm}
\begin{rem}
This assumption enables the value function to be a RCLF.
\end{rem}

In contrast to \cite{Primbs1996nonlinear}, input and disturbance constraints are considered.
The converse constrained optimal control problem for the undisturbed case is discussed in \cite{nevistic1996constrained}.

For Euclidean norm constraints on disturbances and input controls, i.e $\Vert w \Vert \leq \alpha_1$ and $ \Vert u \Vert \leq \alpha_2$, the following optima can be deduced from the existing literature for constrained optimal control problem with disturbances (see \cite{nevistic1996constrained}, Theorem 6).
\begin{equation}
\begin{split}
w^{*} (x) &= \operatorname{sat}_{\alpha_1} \left( \frac{1}{2} g_{1}^{\top}(x) V_{x} \right), \\
u^{*} (x) &= - \operatorname{sat}_{\alpha_2} \left( \frac{1}{2} g_{2}^{\top}(x) V_{x} \right).
\end{split}
\end{equation}
Using these optima in the HJI equation (\ref{eq:HJI}) yields
\begin{equation}
\begin{split}
\left \Vert u^* + g_{2}^{\top}  V_{x} \right \Vert^{2}  - \left \Vert w^* - g_{1}^{\top}  V_{x} \right \Vert^{2} \\
+ V_{x}^\top f+\frac{1}{4} V_{x}^\top \left(g_{1} g_{1}^{\top}-g_{2} g_{2}^{\top}\right) V_{x}+q=0.
\label{eq:HJI_opt}
\end{split}
\end{equation}
Notice that in order to exploit the original converse HJB equation that is presented in \cite{Primbs1996nonlinear}, the two first terms of (\ref{eq:HJI_opt}) must vanish. 
\begin{rem}
From this point onward, state dependencies will be omitted to simplify the notation.
\end{rem}
\begin{asm}
\label{asm:state}
For all $t \geq 0$, $x \in \X$ where: 
\begin{equation}
\begin{split}
\X = \left\{ x \in R^n \mid \left \Vert \frac{1}{2} g_{2}^{\top}  V_{x} \right \Vert < \alpha_2  \cap  
 \left \Vert \frac{1}{2} g_{1}^{ \top} V_{x} \right  \Vert < \alpha_1 \right  \},
\end{split}
\end{equation}
\end{asm}
where $\Vert w \Vert \leq \alpha_1$ and $ \Vert u \Vert \leq \alpha_2$.

This assumption is important because it ensures that the constraints remain inactive, which enables the expression of the optima without the saturation operator. 
Thus, the optima become \cite{Primbs1996nonlinear}:
\begin{equation}
\begin{split}
u^{*} =-\frac{1}{2} g_{2}^{\top} V_{x}, \, \, \, \, 
 &w^{*} =\frac{1}{2} g_{1}^{\top} V_{x},
\end{split}
\label{eq:uoptima}
\end{equation}
which transforms the HJI equation into:
\begin{equation}
V_{x}^\top f+\frac{1}{4} V_{x}^\top \left(g_{1} g_{1}^{\top}-g_{2} g_{2}^{\top}\right) V_{x}+q=0.
\label{eq:HJI_new}
\end{equation}
Now that the framework is set up, the converse optimality problem can be reformulated by adding a robust stabilizability condition.
\paragraph*{Problem statement}
 Given a performance objective $J$, a value function $V$, disturbance coupling functions $g_1$, input coupling functions $g_2$, find a class of \textbf{robustly stabilizable nonlinear systems} $f$ for which this $V$ is the solution of the optimal control problem (\ref{eq:HJI}) over a set $\X$.
 
 The solution to this problem first requires generalizing the existing theory on systems with disturbances, which is presented in the next section.

 \section{Converse Optimality for Systems with Disturbances}
In \cite{Primbs1996nonlinear}, \cite{nevistic1996constrained} and \cite{nevistic1996optimality} the authors assumed $g_1(x)=0$ and presented the converse optimality theory in the undisturbed case  only.
In this section, the admissibility theorem \cite{nevistic1996constrained} and internal dynamics design is extended to systems with disturbances.
\begin{crl}
\label{crl:admissibility}
The pair $(V,q)$ is admissible if and only if: $q=V_x^{\top}h$, for $h \in C$, $h(0)=0$.
\end{crl}

\begin{proof}
This proof has been inspired by the work of Nevisti{\'c} et al. in \cite{nevistic1996optimality}.

\paragraph{Necessity} $(V,q)$ is admissible, then there exists $f$, $g_1$ and $g_2 \in C$ such that:
\begin{equation}
\label{eqn:crl1proof1}
\begin{split}
V_x^{\top}f &+ \frac{1}{4}V_x^{\top} \left(g_1g_1{\top} - g_2g_2^{\top} \right)V_x + q = 0,\\
q &= -V_x^{\top}f + \frac{1}{4}V_x^{\top} \left(g_2g_2^{\top} - g_1g_1^{\top} \right)V_x,\\
q &= V_x^{\top} \left( \frac{1}{4} \left(g_2g_2^{\top} - g_1g_1^{\top}\right)V_x -f \right),\\
q &=V_x^{\top}h,
\end{split}
\end{equation}
where $h = \frac{1}{4}\left(g_2g_2^{\top} - g_1g_1^{\top} \right)V_x - f \in C$ and $h(0) = 0$ from $V_x(0) = 0$ and $f(0)=0$.

\paragraph{Sufficiency} Assume $q=V_x^{\top}h$, with $h \in C$ and $h(0)=0$, for any $g_1$ and $g_2 \in C$:
\begin{equation}
\label{eqn:crl1proof2}
\begin{split}
     V_x^{\top}f + \frac{1}{4}V_x^{\top} \left(g_1g_1{\top} - g_2g_2^{\top}\right)V_x + V_x^{\top}h = 0,\\
     V_x^{\top} \left( f + \frac{1}{4} \left(g_1g_1{\top} - g_2g_2^{\top}  \right) V_x + h \right) = 0,\\
    f + \frac{1}{4} \left(g_1g_1{\top} - g_2g_2^{\top}\right) V_x + h = 0,\\
     f = \frac{1}{4} \left(g_2g_2^{\top} - g_1g_1^{\top}\right) V_x - h.
\end{split}
\end{equation}
$f\in C$ and $f(0) = 0$, therefore $(V,q)$ is admissible.
\end{proof}
\begin{crl}
\label{crl:internaldynamics}
If $(V,q)$ is an admissible pair then the internal dynamics of the converse problem is given by: 
\begin{equation}
         f = \frac{1}{4} \left(g_2g_2^{\top} - g_1g_1^{\top} \right)V_x - h.
         \label{eqn:internal_dynamics}
\end{equation}
\end{crl}
\vspace{0.1cm}
The proof is straightforward from above.
\subsection{Robust Stabilizability}
The main goal behind this work is to ensure that the internal dynamics of the system found in (\ref{eqn:internal_dynamics}) is robustly stabilizable. 
Which is why the RCLF, a generalization of the Lyapunov function for systems with disturbances, is exploited. 
Freeman et al. prove in \cite{FreemanBirkhauser1996} that the existence of a RCLF for any system is a necessary and sufficient condition for robust stabilizability.
This condition is added to the converse optimality problem to achieve robust stabilizability.
Since in converse optimality, the value function is known, it can be chosen such that it is also a RCLF, i.e $C^1$ function in $\mathcal{A}_\infty(\X)$. 
The challenge remaining is to find a proper design for $f$ that satisfies both (\ref{eqn:internal_dynamics}) and (\ref{eq:rclf}).
Therefore, the degree of freedom $h$ in (\ref{eqn:internal_dynamics}) is exploited.

For simplicity, only Euclidean norm constraints on disturbances and control inputs are considered
i.e. the images of the  set valued maps $\mathcal{U}(x)$ and $\mathcal{W}(x)$ are open norm balls and will be noted $\mathcal{U}$ and $\mathcal{W}$ from now on. The optimization problem (\ref{eq:rclf}) is solved in this context as follows. 

Consider the affine system (\ref{eq:sys}) with the performance objective (\ref{eq:perform}) and a value function $V$ that is $C^1$, in $\mathcal{A}_\infty(\X)$ and strongly increasing. 
The inf-sup problem on the Lyapunov derivative of the system is
\begin{equation}
    \inf_{u\in \mathcal{U}} \sup_{w\in \mathcal{W}} L_zV(x,u,w) =  \inf_{u\in \mathcal{U}} \sup_{w\in \mathcal{W}} V_x^\top \left(f + g_1w + g_2u \right),
    \label{eqn:inf_cond}
\end{equation}
where $\mathcal{W} = \{ w \in \mathbb{R}^m \mid \Vert w \Vert < \alpha_1 \}$  and
$\mathcal{U} = \{ u  \in \mathbb{R}^p \mid \Vert u \Vert < \alpha_2 \}$ such that $\alpha_{1,2} > 0$.

Setting up $z=V_x^\top \left(f + g_1w + g_2u \right)$ \red{:}
\begin{equation}
\begin{split}
    \frac{\partial z}{\partial w} =g_1^\top V_x, \, \, \,\,
    &\frac{\partial z}{\partial u} =g_2^\top V_x.    
\end{split}
\end{equation}
Therefore, the optima satisfy 
\begin{equation}
    \begin{split}
        w^\circ = \alpha g_1^\top V_x, \, \, \, \,
        &u^\circ = - \beta g_2^\top V_x,
    \end{split}
\end{equation}
where $\alpha \geq 0$ and $\beta \geq 0$.

Considering that $\Vert \alpha g_1 ^\top V_x \Vert = \alpha_1 $ and  $\Vert \beta g_2 ^\top V_x \Vert = \alpha_2 $,
\begin{equation}
\begin{split}
    w^\circ = \frac{\alpha_1}{\Vert g_1^\top V_x \Vert}g_1^\top V_x, \, \, \, \,
    u^\circ = - \frac{\alpha_2}{\Vert g_2^\top V_x \Vert}g_2^\top V_x.
\end{split}
\end{equation}
Now that the optimization problem is solved for $u$ and $w$, plugging these optima in (\ref{eqn:inf_cond}) yields
\begin{equation}
\begin{split}
  \inf_{u\in \mathcal{U}} \sup_{w\in \mathcal{W}} L_zV(x,u,w) &= V_x^\top f + \frac{\alpha_1}{\Vert g_1^\top V_x \Vert}V_x^\top g_1 g_1^\top V_x \\ &-\frac{\alpha_2}{\Vert g_2^\top V_x \Vert}V_x^\top g_2 g_2^\top V_x. 
  \label{eq:13}
\end{split}
\end{equation}
Introducing the converse optimality's internal dynamics design of (\ref{eqn:internal_dynamics}), it follows that
\begin{equation}
\begin{split}
    &\inf_{u\in \mathcal{U}} \sup_{w\in \mathcal{W}} L_z V(x,u,w) =   \frac{1}{4} V_x^\top(g_2g_2^{\top} - g_1g_1^{\top})V_x - V_x^\top h \\
    & + \frac{\alpha_1}{\Vert g_1^\top V_x \Vert}V_x^\top g_1 g_1^\top V_x - \frac{\alpha_2}{\Vert g_2^\top V_x \Vert}V_x^\top g_2 g_2^\top V_x,
    \end{split}
\end{equation}
and after simplification:
\begin{equation}
\begin{split}
 & \inf_{u\in \mathcal{U}} \sup_{w\in \mathcal{W}} L_z V(x,u,w) =   \frac{1}{4} \Vert g_2^\top V_x \Vert^2 - \alpha_2 \Vert g_2^\top V_x \Vert \\
 &- \frac{1}{4} \Vert g_1^ \top V_x \Vert ^2 + \alpha_1 \Vert g_1 ^ \top V_x \Vert - V_x^\top h.
 \label{eqn:optsimplified}
\end{split}
\end{equation}
According to Proposition \ref{prp:RCLFexistence}, for the class of systems designed from converse optimality to be robustly stabilizable, there must exist a constant $c>0$ such that this double optimization is negative for all $x \in \{ x \in \X \mid V(x) > c\}.$ 
Therefore, a design for $h$ needs to be constructed to satisfy the negativity of the right hand side of (\ref{eqn:optsimplified}) for all $x \in \{ x \in \X \mid V(x) > c\}$.
Before the result is stated, the following lemma is required.
\begin{lem}
\label{lem:valuefunction}
Consider a $C^1$ function $V \in \mathcal{A}_{\infty}(\mathbb{X})$, for every $x \in \mathbb{X}$ that satisfies: $\Vert V_x \Vert > b$ such that $b$ is a strongly positive constant, there exists a constant $c(b) \geq 0$ such $V(x) > c(b)$.
\end{lem}
The proof is straightforward.

\begin{thm}
\label{thm:RSdesign}
Consider a robust converse optimality problem with a value function $V$. 
Let $\mathcal{U}$ and $\mathcal{W}$ be the constraints sets for the inputs and disturbances respectively, and let Assumptions \ref{asm:U}-\ref{asm:state} hold. If the following is satisfied:
\begin{equation}
    h= \frac{1}{4}g_2g_2^\top V_x + \frac{1}{2}g_1g_1^\top V_x + PV_x + \gamma,
    \label{eqn:h_design}
\end{equation}
where $V_x^\top \gamma = 0 $ and 
    \begin{equation}
    P= \alpha_2 g_2 g_2^\top + \frac{\alpha_1 ^2 }{3} \left(1 + \frac{1}{b^2} \right)\frac{1}{1 + V_x^\top V_x}I_{n} + E,
    \label{eqn:design_P}
\end{equation}
where $I_n$ is the identity matrix of dimension $n$, $E$ is a positive definite matrix and $b$ is a positive constant.
Then, the converse system (\ref{eq:sys}) is robustly stabilizable.
\end{thm}
\begin{proof}
To establish this proof,  it is enough to show that the value function is a RCLF.
It follows that the right hand side of (\ref{eqn:optsimplified}) must be negative for all  $x \in \{ x \in \X \mid V(x) > c\}$ where $c$ is a positive constant.

Plugging (\ref{eqn:h_design}) in (\ref{eqn:design_P}) yields
\begin{equation}
\begin{split}
    & \inf_{u\in \mathcal{U}} \sup_{w\in \mathcal{W}} L_z V(x,u,w) =  \frac{1}{4} \Vert g_2^\top V_x \Vert ^2 - \alpha_2 \Vert g_2^\top V_x \Vert \\
    &- \frac{1}{4} \Vert g_1^ \top V_x \Vert ^2 + \alpha_1 \Vert g_1 ^ \top V_x \Vert - \frac{1}{4} V_x^\top g_2g_2^\top V_x \\
    &- \frac{1}{2}V_x^\top g_1g_1^\top V_x - V_x^\top P V_x + V_x^\top \gamma,
\end{split}
\end{equation}
and after modification:
\begin{equation}
\begin{split}
&\inf_{u\in \mathcal{U}} \sup_{w\in \mathcal{W}} L_z V(x,u,w) = - \frac{3}{4} \Vert g_1^\top V_x \Vert ^ 2 + \alpha_1 \Vert g_1 ^ \top V_x \Vert \\ &- \alpha_2 \Vert g_2^\top V_x \Vert - V_x^\top PV_x.
\label{eqn:binomial}
\end{split}
\end{equation}
Notice that (\ref{eqn:binomial}) is a quadratic function of $Y = \Vert g_1^\top V_x \Vert $ that can be written as
\begin{equation}
\label{eqn:F_y}
    F(Y) =  - \frac{3}{4} Y^2 + \alpha_1 Y - \alpha_2 \Vert g_2^\top V_x \Vert - V_x^\top PV_x,
\end{equation}
and it is known that
\begin{equation}
\label{eqn:maxF}
    \max_{Y} \{ F(Y) \} = \frac{\alpha_1 ^2}{3} - \alpha_2\Vert g_2^\top V_x \Vert - V_x^\top P V_x.
\end{equation}
The function (\ref{eqn:F_y}) is a downward open parabola. 
Thus, if the maximum of $F$ (\ref{eqn:maxF}) is negative then the function is negative in $\mathbb{R}$.
Therefore, the construction of $P$ must satisfy
\begin{equation}
    \frac{\alpha_1 ^2}{3} - \alpha_2\Vert g_2^\top V_x \Vert < V_x^\top P V_x, \; \forall x \in \{x \in \X \mid V(x) > c\}.
    \label{eqn:condition_p}
\end{equation}
It remains to show that (\ref{eqn:design_P}) satisfies this condition.
Plugging (\ref{eqn:design_P}) in the right hand side of (\ref{eqn:condition_p}) gives
\begin{equation}
\begin{split}
    V_x^\top P V_x &= V_x^\top \alpha_2 g_2 g_2^\top V_x \\ &+ V_x^\top \frac{\alpha_1 ^2 }{3}\left(1 + \frac{1}{b^2} \right)\frac{1}{1 + V_x^\top V_x}I_{n} V_x  + V_x^\top E V_x,
\end{split}
\end{equation}
which can be written as
\begin{equation}
\begin{split}
   V_x^\top P V_x  &= \alpha_2 \Vert g_2^\top V_x \Vert^2 + \frac{\alpha_1 ^2}{3}\left(1 + \frac{1}{b^2}\right)\frac{\Vert V_x\Vert ^2}{1 + \Vert V_x \Vert ^2} \\
&+ \Vert V_x^\top E V_x\Vert.
\end{split}
\end{equation}
On the other hand,
\begin{equation}
 - \alpha_2\Vert g_2^\top V_x \Vert < \alpha_2 \Vert g_2^\top V_x \Vert^2, \, \, \forall x \in \mathbb{X}.
 \label{eqn:first}
\end{equation}

Let $b > 0$ such that : $\Vert V_x \Vert > b$, then, according to Lemma \ref{lem:valuefunction}, there exists a constant $c(b)  \geq 0$ such that : $V(x) > c(b)$.

As per Assumption \ref{asm:valuefunc}:
\begin{equation}
        \Vert V_x \Vert > b,
\end{equation}
inverting and squaring both sides yields
\begin{equation}
        \frac{1}{\Vert V_x \Vert^2} < \frac{1}{b^2}
\end{equation}
adding 1 to each side yields
\begin{equation}
        \frac{\Vert V_x \Vert^2 +1 }{\Vert V_x \Vert^2} < 1 + \frac{1}{b^2},
\end{equation}
multiplying by the positive constant $\frac{\alpha_1 ^2}{3}$:
\begin{equation}
        \frac{\alpha_1^2}{3}\left(\frac{\Vert V_x \Vert^2 +1 }{\Vert V_x \Vert^2}\right) < \frac{\alpha_1^2}{3}\left(1 + \frac{1}{b^2} \right),
\end{equation}
rearranging gives
\begin{equation}
        \frac{\alpha_1^2}{3}<\frac{\alpha_1^2}{3}\left(1 + \frac{1}{b^2}\right)\frac{\Vert V_x \Vert^2}{1 + \Vert V_x \Vert^2},
\end{equation}
for all $x \in \{ x \in \X \mid V(x) > c(b) \}$.

Adding (\ref{eqn:first}) to the latter,
\begin{equation}
\begin{split}
   - \alpha_2\Vert g_2^\top V_x \Vert + \frac{\alpha_1^2}{3}<\frac{\alpha_1^2}{3}\left(1 + \frac{1}{b^2}\right)\frac{\Vert V_x\Vert^2}{1 + \Vert V_x\Vert^2} \\+ \alpha_2 \Vert g_2^\top V_x \Vert^2,
\end{split}
\end{equation}
and since $E$ is a positive definite matrix:
\begin{equation}
\begin{split}
    &- \alpha_2\Vert g_2^\top V_x \Vert + \frac{\alpha_1^2}{3} < \frac{\alpha_1^2}{3}\left(1 + \frac{1}{b^2}\right)\frac{\Vert V_x \Vert^2}{1 + \Vert V_x\Vert^2} \\&+ \alpha_2 \Vert g_2^\top V_x \Vert^2 +  \Vert V_x^\top E V_x\Vert,  \; \forall x \in \{ x \in \X \mid V(x) > c(b) \},
\end{split}
\end{equation}
which is equivalent to (\ref{eqn:condition_p}).
\end{proof}

This theorem allows the construction of robustly stabilizable nonlinear internal dynamics in the context of converse optimality.
The stability analysis in open loop of the resulting class of system is provided in the next subsection.
\subsection{Stability in Open Loop}
In this subsection, the proof of the asymptotic open loop stability of the robustly stabilizable class of nonlinear systems obtained from converse optimality is presented.
This is shown by proving that the value function is also a Lyapunov function for the system.

The open loop system can be written this way:
\begin{equation}
    \dot{x}=f(x), \, \, x \in \X,
\end{equation}
where, with the combination of (\ref{eqn:internal_dynamics}), (\ref{eqn:h_design}) and (\ref{eqn:design_P}), the internal dynamics of the system is
\begin{equation}
\begin{split}
f &= - \frac{3}{4}g_1 g_1^\top V_x - \alpha_2 g_2g_2^\top V_x \\ &- \frac{\alpha_1 ^2}{3}\left(1 + \frac{1}{b^2}\right)\frac{1}{1+ V_x^\top V_x}V_x - EV_x + \gamma.
\end{split}
\label{eq:fexp}
\end{equation}
The value function $V$ is positive $\forall x \in \X \setminus\{0\}$ by definition, which makes it a candidate Lyapunov function.
Consider its Lie derivative
\begin{equation}
    \dot{V}=V_x^\top f,
\end{equation}
replacing $f$ by its expression (\ref{eq:fexp}):
\begin{equation}
\begin{split}
    \dot{V} &=  - \frac{3}{4}V_x^\top g_1 g_1^\top V_x - \alpha_2 V_x^\top g_2g_2^\top V_x \\ &- \frac{\alpha_1 ^2}{3}\left(1 + \frac{1}{b^2}\right)\frac{1}{1 + V_x^\top V_x}V_x^\top V_x - V_x^\top EV_x + V_x^\top \gamma, 
    \end{split}
\end{equation}
where $V_x^\top \gamma = 0$ by definition, and $E$ is a positive definite matrix. 
Therefore the terms $ \frac{\alpha_1 ^2}{3}\left(1 + \frac{1}{b^2}\right)\frac{1}{1 + V_x^\top V_x}V_x^\top V_x - V_x^\top EV_x$ are always negative for all $x \in \mathbb{X} \setminus \{0\}$ then, it is straightforward that $\dot{V} < 0$ for all $x \in \mathbb{X} \setminus \{0\}$.
Thus, and according to Lyapunov theorem, the designed systems are asymptotically stable \cite{khalil2002nonlinear}.

\section{Case Study}

\subsection{ $2 \times 2$ System with Constant Coupling Functions}
The motivation of this example is to give a generalized class of systems designed by this approach. 
For simplicity and to avoid complex calculations, 2D systems are chosen. The coupling functions are unknown and assumed constant. 
In principle, this can also be achieved for higher system degrees and state varying coupling functions.

Consider $g_1 = \begin{bmatrix}g_{11} & g_{12} \end{bmatrix} ^\top $, $g_2 = \begin{bmatrix}g_{21} & g_{22} \end{bmatrix} ^\top $ where $g_{ij} \in \R$ for $i,j \in \{1 , 2\}$. $E=\begin{bmatrix}e_1 & 0 \\ 0 & e_2 \end{bmatrix}$ where $e_{1,2} \in \R_{\geq 0}$,
and the quadratic value function $V=\frac{1}{2}x_1^2 + \frac{1}{2}x_2^2$ and $\gamma= \begin{bmatrix}-x_2 & x_1\end{bmatrix}^\top$ .

Using Theorem \ref{thm:RSdesign}, the following system is obtained
\begin{equation}
\begin{aligned}
\label{eqn:2Dsys}
  f_1(x) =\\
   &\left[ \begin{matrix} -\frac{3}{4}g_{11}^2 - \alpha_2g_{21}^2 - \frac{k}{1 + x_1^2 + x_2^2}- e_1 \\ -\frac{3}{4}g_{12}g_{11} - \alpha_2g_{22}g_{21} + 1   \end{matrix}\right.\\
&\qquad \qquad
\left.\begin{matrix}
     {}-\frac{3}{4}g_{11}g_{12} - \alpha_2 g_{21}g_{22} -1
     \\ - \frac{3}{4}g_{12}^2 - \alpha_2 g_{22}^2 - e_2 - \frac{c}{1 + x_1^2 + x_2^2}
    \end{matrix} \right]x,
\end{aligned}
\end{equation}
where $k = \frac{\alpha_1^2}{3}\left(1 + \frac{1}{b^2}\right)$.

\begin{figure}[h]
\centering
    \includegraphics[scale=0.41]{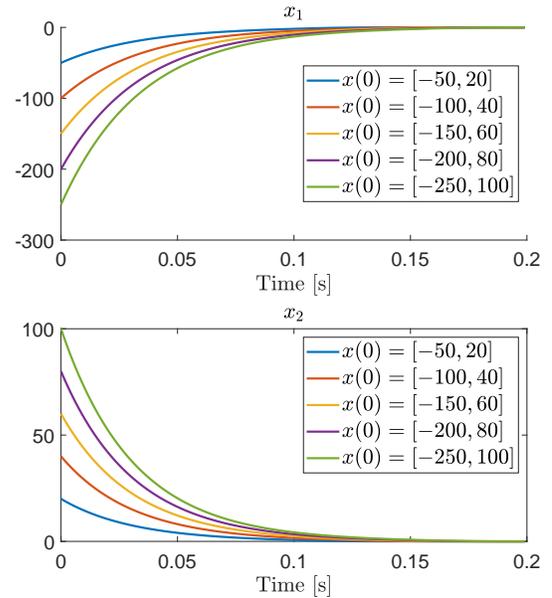}
    \caption{Transient behaviour of (\ref{eqn:2Dsys}) with varying initial condition}
    \label{fig:2D}
\end{figure}

For numerical illustration, Fig. \ref{fig:2D} shows the the open loop response of this design for different initial conditions, with the following parameters:
$g_1 = \begin{bmatrix}1 & 5 \end{bmatrix} ^\top $,  $g_2 = \begin{bmatrix}-1 & 0 \end{bmatrix} ^\top $, $E=\begin{bmatrix}10 & 0 \\ 0 & 20 \end{bmatrix}$, $\alpha_1 = 10$, $\alpha_2 = 20$, $b=\sqrt{0.5}$.
\vspace{-1em}

The asymptotic stability of the open loop system can also be observed as both states asymptotically converge to the origin.

\subsection{ $3 \times 3$ System with Nonlinear Coupling Functions}
To get a highly nonlinear system, this example considers the following nonlinear coupling functions: $g_1= \begin{bmatrix} -x_2^2 & x_1x_2 & x_3 \end{bmatrix}^\top$ and $g_2 = \begin{bmatrix} x_3 & 1 & -x_2 \end{bmatrix}^\top$ and the following parameters: $E = \begin{bmatrix}10 & 0 & 0 \\ 0 & 5 & 0 \\ 0 & 0 & 5  \end{bmatrix}$, $b=\sqrt{0.5}$, $\alpha_1 = 5$, $\alpha_2 = 10$, keeping a quadratic value function $V=\frac{1}{2}x_1^2 + \frac{1}{2}x_2^2 + \frac{1}{2}x_3^2$, $\gamma = \begin{bmatrix}-x_2 & x_1 & 0 \end{bmatrix}^\top$ and $x_0 = \begin{bmatrix} 5 & 4 & -1 \end{bmatrix}^\top$.
The resulting internal dynamics is
\begin{equation}
\begin{aligned}
 & f_2(x) =\\
   &\left[ \begin{matrix} -\frac{3}{4}x_1x_2^4 + \frac{3}{4}x_1x_2^4 + \frac{3}{4}x_2^2x_3^2 - 10x_1x_3^2 \\
   \frac{3}{4}x_1^2x_2^3 - \frac{3}{4}x_1^2x_2^3 - \frac{3}{4}x_1x_2x_3^2 - 10x_1x_3\\
   \frac{3}{4}x_1x_2^2x_3 - \frac{3}{4}x_1x_2^2x_3 - \frac{3}{4}x_3^3 + 10x_1x_2x_3 + 10x_2^2
   \end{matrix}\right.\\
&\qquad
\left.\begin{matrix}
-10x_2x_3 + 10x_2x_3^2 - \frac{25x_1}{1 + x_1^2 + x_2^2 + x_3^2}-10x_1 - x2\\
+ 10x_2x_3 - \frac{25x_2}{1 + x_1^2 + x_2^2 + x_3^2} - 15x_2 + x_1\\
- 10x_2^2x_3- \frac{25x_3}{1 + x_1^2 + x_2^2 + x_3^2} - 5x_3
    
\end{matrix} \right].
\end{aligned}
\end{equation}
The matrix $E$ can also be chosen as state-dependent as long as it remains positive definite. Such a choice would increase even more the nonlinear aspect of the internal dynamics $f_2$.
\vspace{-1em}
\begin{figure}[h]
 \centering
    \includegraphics[scale=0.5]{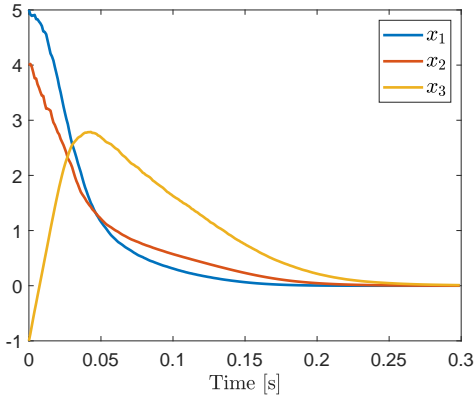}
    \caption{State convergences of a 3D system}
    \label{fig:3Dx}
\end{figure}
\vspace{-1em}
\begin{figure}[h]
\centering
    \includegraphics[scale=0.40]{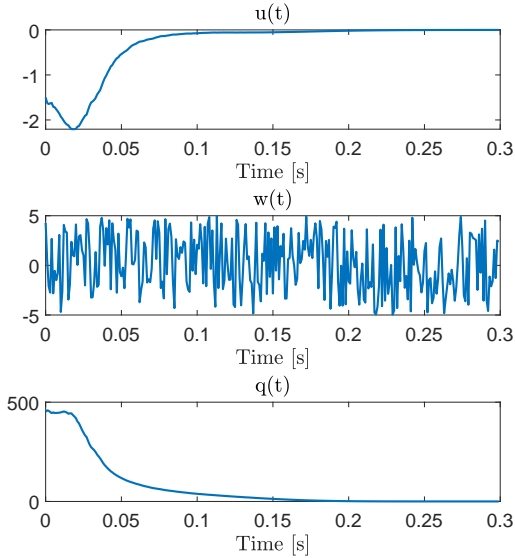}
    \caption{Control, disturbance and utility functions for a 3D system}
    \label{fig:3Duwq}
\end{figure}
\vspace{-1em}

In Fig. \ref{fig:3Dx} the previous system is simulated in closed loop, controlled with $u^{*}$ in (\ref{eq:uoptima}) with a uniformly distributed disturbance in $\begin{bmatrix}-5 & 5\end{bmatrix}$. 
The asymptotic convergence to the origin can be observed after a lightly disturbed transient phase. 
Even though the amplitude of the disturbances is significant, it is smoothed and barely affects the states trajectory and control input's variation. 
The latter is represented in Fig. \ref{fig:3Duwq} along with the disturbance and utility functions.

\section{Conclusion}
This work extends the converse optimality problem to systems with disturbances. Previous converse optimality results are generalized for this framework. 
A robustly stabilizable class of nonlinear systems is designed utilizing a robust control Lyapunov function.
Combining, thus, converse optimality and inverse optimality theories locally.

Future work will focus on stabilizing practical systems close to robust conversely designed systems.
Additionally, data driven approaches will be investigated for the case where only measurements are available. 

\addtolength{\textheight}{-10cm}   





\phantomsection
\addcontentsline{toc}{chapter}{Bibliography} 
\bibliographystyle{acm}
\bibliography{bib}

\end{document}